\documentclass[11pt,fleqn]{article}
\usepackage{amsmath,amssymb,graphicx}

\begin{document}
\sloppy
\newtheorem{axiom}{Axiom}[section]
\newtheorem{conjecture}[axiom]{Conjecture}
\newtheorem{corollary}[axiom]{Corollary}
\newtheorem{definition}[axiom]{Definition}
\newtheorem{example}[axiom]{Example}
\newtheorem{fact}[axiom]{Fact}
\newtheorem{lemma}[axiom]{Lemma}
\newtheorem{observation}[axiom]{Observation}
\newtheorem{open}[axiom]{Problem}
\newtheorem{proposition}[axiom]{Proposition}
\newtheorem{theorem}[axiom]{Theorem}

\renewcommand{\topfraction}{1.0}
\renewcommand{\bottomfraction}{1.0}

\newcommand{\proof}{\emph{Proof.}\ \ }
\newcommand{\qed}{~~$\Box$}
\newcommand{\rz}{{\mathbb{R}}}
\newcommand{\nz}{{\mathbb{N}}}
\newcommand{\zz}{{\mathbb{Z}}}
\newcommand{\eps}{\varepsilon}
\newcommand{\cei}[1]{\lceil #1\rceil}
\newcommand{\flo}[1]{\left\lfloor #1\right\rfloor}

\newcommand{\aaa}{\alpha}
\newcommand{\bbb}{\beta}
\newcommand{\ccc}{\gamma}

\newcommand{\lap}{\mbox{LAP}}
\newcommand{\qap}{\mbox{QAP}}
\newcommand{\xxx}[1]{^{\scriptscriptstyle[#1]}}

\title{{\bf Linearizable special cases of the QAP}}
\author{
\sc Eranda \c{C}ela\thanks{{\tt cela@opt.math.tu-graz.ac.at}.
Institut f\"ur Optimierung und Diskrete Mathematik, TU Graz, Steyrergasse 30, A-8010 Graz, Austria (corresponding author)}
\and\sc Vladimir G.\ Deineko\thanks{{\tt Vladimir.Deineko@wbs.ac.uk}.
Warwick Business School, The University of Warwick, Coventry CV4 7AL, United Kingdom}
\and\sc Gerhard J.\ Woeginger\thanks{{\tt gwoegi@win.tue.nl}.
Department of Mathematics and Computer Science, TU Eindhoven, P.O.\ Box 513,
5600 MB Eindhoven, Netherlands}
}
\date{}
\maketitle

\begin{abstract}
We consider special cases of the quadratic assignment problem (QAP) that are linearizable in the 
sense of Bookhold.
We provide combinatorial characterizations of the linearizable instances of the weighted feedback 
arc set QAP, and of the linearizable instances of the traveling salesman QAP. 
As a by-product, this yields a new well-solvable special case of the weighted feedback arc set problem. 

\end{abstract}

\medskip
\section{Introduction}
The \emph{Quadratic Assignment Problem} (QAP) and the \emph{Linear Assignment Problem} (LAP) are 
two important and well-studied problems in combinatorial optimization; we refer the reader to the 
books by \c{C}ela \cite{Cela-book} and by Burkard, Dell'Amico \& Martello \cite{Burkard-book} for 
comprehensive surveys on these problems.
The QAP in Koopmans-Beckmann form \cite{KoBe1957} takes as input two $n\times n$ square matrices
$A=(a_{ij})$ and $B=(b_{ij})$ with real entries, and assigns to every permutation $\pi\in S_n$ 
(where $S_n$ denotes the set of permutations of $\{1,2,\ldots,n\}$) the corresponding objective value
\begin{equation}
\label{eq:qap}
\qap(A,B,\pi) ~:=~ \sum_{i=1}^n\sum_{j=1}^n~ a_{\pi(i)\pi(j)} \, b_{ij}.
\end{equation}
The LAP takes as input a single $n\times n$ matrix $C=(c_{ij})$, and assigns to every 
permutation $\pi\in S_n$ the objective value
\begin{equation}
\label{eq:lap}
\lap(C,\pi) ~:=~ \sum_{i=1}^n~ c_{i\pi(i)}.
\end{equation}
The usual goal in these optimization problems is to identify permutations $\pi$ that minimize 
the objective values (\ref{eq:qap}) and (\ref{eq:lap}), respectively.
The QAP is NP-hard and extremely difficult to solve, whereas the LAP is polynomially solvable 
and fairly harmless \cite{Cela-book,Burkard-book}.

Bookhold \cite{Bo1990} calls an instance of the QAP (that is, two $n\times n$ matrices $A$ and $B$) 
\emph{linearizable}, if there exists a corresponding instance of the LAP (that is, a single 
$n\times n$ matrix $C$) such that
\begin{equation}
\label{eq:linearizable}
\qap(A,B,\pi) ~=~ \lap(C,\pi) \mbox{\qquad for all permutations $\pi\in S_n$.}
\end{equation}
Of course linearizable instances of the QAP are polynomially solvable by simply solving the 
corresponding instance of the LAP.

In a tour de force, Kabadi \& Punnen \cite{KaPu2011,PuKa2013} designed an $O(n^2)$ polynomial
time algorithm for recognizing linearizable instances of the QAP in Koopmans-Beckmann form.
Furthermore, \cite{PuKa2013} derived a purely combinatorial characterization of all linearizable
QAP instances with \emph{symmetric} matrices $A$ and $B$: such instances are linearizable if 
and only if one of the two matrices is a weak sum matrix (see Section~\ref{sec:qap} for a more 
precise statement of this result).
Hence linearizable \emph{symmetric} QAP instances are fully understood and carry a highly
restrictive combinatorial structure.
The structure of \emph{asymmetric} linearizable QAP instances is much richer, and it seems to 
be very difficult to extend the algorithmic characterization of \cite{PuKa2013} to a clean 
combinatorial characterization.
Asymmetric linearizable QAPs are the topic of the present paper.

\paragraph{Results of this paper.}
We perform a combinatorial study on Bookhold linearizations of two prominent and well-studied 
families of \emph{asymmetric} QAP instances: the feedback arc set problem (FAS) and the traveling 
salesman problem (TSP). 
As our main results, we derive the following combinatorial characterizations for these problems.
\begin{itemize}
\item
An instance of the FAS is linearizable if and only if in the underlying arc weight matrix all 
the 3-cycles are balanced; this means that for every cycle on three vertices, the total weight 
of its clockwise traversal equals the total weight of its counter-clockwise traversal.
\item
An instance of the TSP is linearizable if and only if the underlying distance matrix is a weak 
sum matrix; this means that the (asymmetric) distances from city $i$ to city $j$ are given as the 
sum of two parameters, one of which only depends on $i$ while the other one only depends on $j$.
\end{itemize}

For the TSP, our results indicate that linearizations will not lead to new well-solvable instances.
In fact linearizations will not be able to add anything new to the TSP literature, as TSP 
instances on weak sum matrices have been fully analyzed a long time ago.
It is known that for weak sum matrices, all feasible solutions yield the same TSP objective value.
Gabovich \cite{Gabovich1976} further showed that weak sum matrices are the \emph{only} matrices 
with that property.

For the FAS, our results indicate that linearizations are sometimes useful.
There is one branch of research on the QAP that concentrates on the algorithmic behavior of 
strongly structured special cases; see for instance Burkard \& al \cite{BCRW1998}, 
Deineko \& Woeginger \cite{DeWo1998}, or \c{C}ela, Deineko \& Woeginger \cite{Cela2012} for 
typical results in this direction.
Our results contribute a new well-solvable case to this research branch.
Our proof method analyzes certain linear combinations of certain simple 0-1 matrices, and
hence is similar in spirit to the approaches in \cite{BCRW1998,DeWo1998,Cela2012}.

\paragraph{Organization of the paper.}
Section~\ref{sec:matrix} summarizes the relevant matrix classes and provides a characterization
of balanced 3-cycle matrices. 
Section~\ref{sec:qap} states several observations and results on linearizable QAPs.
Section~\ref{sec:fas-qap} derives our results on the feedback arc set QAP, and 
Section~\ref{sec:tsp-qap} gives the results on the traveling salesman QAP.
Section~\ref{sec:conclusion} completes the paper with a short conclusion.

\medskip
\section{The central matrix classes}
\label{sec:matrix}
In this section we summarize definitions and results around several matrix classes that will
play a central role in our investigations.
All matrices in this paper have real entries, and most of them are square matrices.
An $n\times n$ matrix $A=(a_{ij})$ is a \emph{sum matrix}, if there exist real 
numbers $\aaa_1,\ldots,\aaa_n$ and $\bbb_1,\ldots,\bbb_n$ such that 
\begin{equation}
\label{eq:sum}
a_{ij}=\aaa_i+\bbb_j \mbox{\qquad for $1\le i,j\le n$.}
\end{equation}
Matrix $A$ is a \emph{weak sum matrix}, if $A$ can be turned into a sum matrix by appropriately
changing the entries on its main diagonal.
Matrix $A$ is a \emph{directed cut matrix}, if there exists a subset $I\subseteq\{1,2,\ldots,n\}$
such that 
\begin{equation}
\label{eq:dicut}
a_{ij} ~=~ \left\{
\begin{array}{cl}
1 &\mbox{\quad if $i\in I$ and $j\notin I$} \\[0.5ex]
0 &\mbox{\quad otherwise}
\end{array} \right.
\end{equation}
In graph theoretic terms, the entries in (\ref{eq:dicut}) encode the arcs of the directed 
cut from vertex set $I$ to the complement of $I$.
We will sometimes say that the directed cut matrix is \emph{induced} by $I$.

Three indices $i,j,k\in\{1,2,\ldots,n\}$ are said to form a \emph{balanced 3-cycle} in 
an $n\times n$ matrix $A$, if the corresponding entries satisfy
\begin{equation}
\label{eq:b3c}
a_{ij}+a_{jk}+a_{ki} ~=~ a_{ji}+a_{kj}+a_{ik}.
\end{equation}
This means that the total weight on the clockwise cycle $i,j,k$ equals the total weight on 
the counter-clockwise cycle $k,j,i$.
Matrix $A$ is a \emph{balanced 3-cycle matrix}, if every three indices $i,j,k$ satisfy (\ref{eq:b3c}).

Note that (\ref{eq:b3c}) trivially holds whenever two of the indices $i,j,k$ coincide. 
Note furthermore that the condition (\ref{eq:b3c}) is linear.
Hence the class of balanced 3-cycle matrices is closed under addition and under multiplication by a scalar, 
and forms a subspace of the space of $n\times n$ matrices.
The following theorem derives a characterization of balanced 3-cycle matrices that is
crucial for our arguments in Section~\ref{sec:fas-qap}.

\begin{theorem}
\label{th:subspace}
An $n\times n$ matrix $A$ is a balanced 3-cycle matrix, if and only if it can be written as 
the sum of a symmetric matrix and a linear combination of directed cut matrices.
\end{theorem}
\proof
For the if part, first observe that any symmetric matrix $A$ trivially satisfies (\ref{eq:b3c}).
Next consider the case of a directed cut matrix $A$ that is induced by $I\subseteq\{1,2,\ldots,n\}$, 
and let $i,j,k$ be three indices.
If all three of $i,j,k$ are contained in $I$ or if none of them is contained in $I$,
then the values of the left hand side and right hand side in (\ref{eq:b3c}) both are~$0$.
If exactly one or two of $i,j,k$ are contained in $I$, then the values of the left hand side 
and right hand side in (\ref{eq:b3c}) both are~$1$.
Hence any symmetric matrix and any directed cut matrix is a balanced 3-cycle matrix, and 
the linearity of (\ref{eq:b3c}) completes the first part of the proof.

For the only if part, we first subtract an appropriately chosen symmetric matrix from matrix $A$
such that afterwards all entries in $A$ are non-negative and satisfy
\begin{equation}
\label{eq:sub.1}
a_{ij}\,a_{ji}=0 \mbox{\qquad for all $i,j$ with $1\le i,j\le n$.}
\end{equation}
We fix two indices $r$ and $s$ such that the value $a_{rs}$ is maximum among all the
entries in matrix $A$. 
If $a_{rs}=0$, then $A$ is the all zero matrix and we are done.
Otherwise $a_{rs}$ is positive, and (\ref{eq:sub.1}) implies $a_{sr}=0$.
We define set $I$ to contain all indices $i$ satisfying $a_{ri}<\frac12 a_{rs}$; 
note that $r\in I$ and $s\notin I$.

Now consider two arbitrary indices $i\in I$ and $j\notin I$, which by definition
fulfill $a_{ri}<\frac12 a_{rs}\le a_{rj}$.
By (\ref{eq:sub.1}) we then have $a_{jr}=0$. 
In case also $a_{ij}=0$ holds, (\ref{eq:b3c}) would yield 
\begin{equation}
\label{eq:sub.2}
\frac12 a_{rs} ~\le~ a_{rj} ~\le~ a_{ji}+a_{rj}+a_{ir} 
~=~ a_{ij}+a_{jr}+a_{ri} ~=~ a_{ri} ~<~ \frac12 a_{rs}.
\end{equation}
This contradiction implies that
\begin{equation}
\label{eq:sub.3}
a_{ij}>0 \mbox{\qquad whenever $i\in I$ and $j\notin I$.}
\end{equation}
Let $A'$ be the directed cut matrix induced by $I$, and let $p\in I$ and $q\notin I$ be 
the indices with the smallest value $a_{pq}$; then $a_{pq}>0$ by (\ref{eq:sub.3}).
The matrix $A-a_{pq}A'$ has non-negative entries, satisfies (\ref{eq:sub.1}), 
and has at least one more zero entry than matrix $A$ (as it also has a zero entry at 
the crossing of row $p$ and column $q$).  

We iterate this step and repeatedly subtract such matrices $a_{pq}A'$ from $A$ and
thereby increase the number of zero entries.
When we finally reach the all zero matrix, the subtracted matrices yield the desired representation 
of $A$ as sum of a symmetric matrix and a linear combination of directed cut matrices.
\qed

\medskip
\section{Linearizations of the QAP}
\label{sec:qap}
In this section we collect some observations and results around linearizable QAPs.
The following statement belongs to the QAP folklore and has been known (in slightly different
formulations) for decades.
\begin{proposition}
\label{pr:sum-qap}
(Folklore)~
If one of the matrices $A$ and $B$ is a weak sum matrix, then the QAP for $A$ and $B$ is linearizable. 
\qed
\end{proposition}

If matrix $A$ in some QAP instance is symmetric, then matrix $B$ may also be made symmetric by
replacing it by $\frac12(B+B^T)$.
Therefore the QAP literature only considers symmetric QAPs (where both matrices are symmetric)
and asymmetric QAPs (where both matrices are asymmetric).
The following result establishes the reverse of Proposition~\ref{pr:sum-qap} for the case of 
symmetric matrices.
\begin{proposition}
\label{pr:KabPun}
(Punnen \& Kabadi \cite{PuKa2013})~
If the QAP for two symmetric matrices $A$ and $B$ is linearizable, then one of $A$ and $B$ 
is a weak sum matrix.
\qed
\end{proposition}

Propositions~\ref{pr:sum-qap} and~\ref{pr:KabPun} provide a full combinatorial characterization 
of linearizable symmetric QAPs.
In strong contrast to this, the structure of asymmetric linearizable QAPs is much richer, and in
particular is not tied to weak sum matrices.
For an illustration, consider the following three matrices:
\begin{equation}
\label{eq:example}
A= \left(
\begin{array}{cccc}
0 &0 &1 &1 \\
0 &0 &1 &1 \\
0 &0 &0 &0 \\
0 &0 &0 &0 
\end{array}\right)
\quad
B= \left(
\begin{array}{cccc}
0 &1 &1 &1 \\
0 &0 &1 &1 \\
0 &0 &0 &1 \\
0 &0 &0 &0 
\end{array}\right)
\quad
C= \left(
\begin{array}{cccc}
0 &1 &2 &3 \\
-1&0 &1 &2 \\
0 &0 &0 &0 \\
0 &0 &0 &0 
\end{array}\right)
\end{equation}
Note that matrices $A$ and $B$ are asymmetric and that neither of them is a weak sum matrix.
Lemma~\ref{le:fas.2} in Section~\ref{sec:fas-qap} yields that the QAP for $A$ and $B$ is
linearizable, and that matrix $C$ is one possible linearization for it.

We close this section with a simple but useful observation.
\begin{lemma}
\label{le:linear}
Let $A_1$, $A_2$ and $B$ be $n\times n$ matrices such that the QAP with matrices
$A_1$ and $B$ as well as the QAP with matrices $A_2$ and $B$ are linearizable.
Then for any real numbers $\lambda_1$ and $\lambda_2$, also the QAP with matrices
$\lambda_1A_1+\lambda_2A_2$ and $B$ is linearizable.
\end{lemma}
\proof
For $k\in\{1,2\}$ let $C_k$ be a matrix such that $\qap(A_k,B,\pi)=\lap(C_k,\pi)$
for all permutations $\pi\in S_n$.
Then $\qap(\lambda_1A_1+\lambda_2A_2,B,\pi)=\lap(\lambda_1C_1+\lambda_2C_2,\pi)$ for all $\pi\in S_n$.
\qed

\medskip
\section{The feedback arc set QAP}
\label{sec:fas-qap}
A \emph{feedback arc set} in a directed graph $G=(V,E)$ is a subset $E'$ of the arcs such that
the subgraph $(V,E-E')$ is a directed acyclic graph; in other words, the subset $E'$ contains 
at least one arc from every directed cycle in $G$.
The goal is to find a feedback arc set of minimum cardinality.
We refer the reader to the survey article \cite{FePaRe2000} by Festa, Pardalos \& Resende
for more information on this problem.

The problem of finding a feedback arc set of minimum cardinality can be modeled as a QAP of 
size $n=|V|$. 
Matrix $A$ is the adjacency matrix of $G$ (so that $a_{ij}=1$ whenever there is an arc from 
vertex $i$ to vertex $j$, and $a_{ij}=0$ otherwise), and matrix $B$ is the $n\times n$ 
\emph{feedback arc} matrix $F_n=(f_{ij})$ whose entries are defined as follows:
\begin{equation}
\label{eq:fas}
f_{ij} ~=~ \left\{
\begin{array}{cl}
1 &\mbox{\quad if~ $1\le j  < i\le n$} \\[0.5ex]
0 &\mbox{\quad if~ $1\le i\le j\le n$}
\end{array} \right.
\end{equation}
In graph theoretic terms, matrix $F_n$ is the adjacency matrix of the directed graph whose
vertices are laid out on the integers $1,2,\ldots,n$, and whose arc set contains all possible
backward arcs (that is, arcs going back from a vertex to another vertex with lower number).
The permutation $\pi$ in the QAP then specifies a topological ordering of the acyclic subgraph 
$(V,E-E')$.
In the corresponding objective value (\ref{eq:qap}), all the forward arcs (from vertices with low
number to vertices with high number) are matched with a $0$ entry in $F_n$ and all the backward arcs
(from vertices with high number to vertices with low number) are matched with a $1$ entry in $F_n$.
The backward arcs form a feedback arc set, and minimizing the cardinality of this set exactly
corresponds to minimizing the objective value of the QAP.

The general \emph{feedback arc set QAP} considers the arc-weighted version, where the goal is 
to find a feedback arc set of minimum weight. 
The first matrix $A$ in the QAP has arbitrary real entries and encodes the arc weights, while 
the second matrix is the feedback arc matrix $F_n$ as specified in (\ref{eq:fas}).
We will call this problem the \emph{FAS-QAP for matrix~$A$}, or just \emph{FAS-QAP} for short.  
The FAS-QAP is NP-hard, as it models the NP-hard feedback arc set problem in directed 
graphs \cite{GaJo1979}.
In the following, we will concisely characterize all linearizable instances of the FAS-QAP.

\begin{lemma}
\label{le:fas.1}
For any symmetric matrix $A$, the FAS-QAP for matrix $A$ is linearizable.
\end{lemma}
\proof
No matter whether vertex $i$ comes before vertex $j$ or after vertex $j$ in the layout, 
the contribution of this vertex pair to the objective function exactly equals $a_{ij}$.
Hence all permutations yield exactly the same objective value for this QAP instance, and 
the instance can be linearized trivially by a matrix $C$ that yields the same constant 
LAP objective value for all permutations.
\qed 

\begin{lemma}
\label{le:fas.2}
For any directed cut matrix $A$, the FAS-QAP for matrix $A$ is linearizable.
\end{lemma}
\proof
We assume without loss of generality that the $n\times n$ directed cut matrix $A$ is
induced by $I=\{1,\ldots,k\}$.
For discussing the FAS-QAP, it is convenient to use the graph theoretic interpretation
described at the beginning of this section.
Consider a permutation $\pi$ that assigns the $k$ vertices of $I$ to the $k$ positions
$p_1<p_2<\cdots<p_k$ in the layout.
Then the vertex assigned to position $p_i$ (with $1\le i\le k$) contributes $p_i-i$ 
backward arcs to the objective value.
Indeed, there are $p_i-1$ positions to the left of position $p_i$, of which $i-1$ are occupied 
by vertices in $I$ while the remaining $p_i-i$ positions are occupied by vertices not in $I$.
There is a backward arc from the vertex at position $p_i$ to each of these $p_i-i$ vertices 
not in $I$.
For the objective value in (\ref{eq:qap}) this yields
\begin{equation}
\label{eq:dcut.1}
\qap(A,F_n,\pi) ~=~ \sum_{i=1}^k (p_i-i).
\end{equation}
For the linearization we use the $n\times n$ matrix $C$ whose first $k$ rows are given by
$c_{ij}=j-i$ for $i=1,\ldots,k$ and $j=1,\ldots,n$, and whose remaining $n-k$ rows only
contain zeroes; see (\ref{eq:example}) for an example.
The objective value in (\ref{eq:lap}) then becomes
\begin{equation}
\label{eq:dcut.2}
\lap(C,\pi) ~=~ \sum_{i=1}^n~ c_{i\pi(i)} ~=~ \sum_{i=1}^k (\pi(i)-i).
\end{equation}
Since the positions $p_1<p_2<\cdots<p_k$ are the values $\pi(1),\ldots,\pi(k)$ ordered by size,
the objective values in (\ref{eq:dcut.1}) and (\ref{eq:dcut.2}) coincide.
\qed

\bigskip
For an $n\times n$ matrix $A$ and a subset $J\subseteq\{1,\ldots,n\}$, the principal submatrix
$A[J]$ results by removing from $A$ all the rows and columns whose index is not in $J$.
\begin{lemma}
\label{le:fas.3}
If the FAS-QAP for an $n\times n$ matrix $A$ is linearizable, then for any 
$J\subseteq\{1,\ldots,n\}$ the FAS-QAP for the principal submatrix $A[J]$ is also linearizable.
\end{lemma}
\proof
We assume without loss of generality that $J=\{1,\ldots,k\}$.
For a permutation $\pi\in S_k$ we define its extension $\pi^+\in S_n$ by $\pi^+(i)=\pi(i)$ for 
$1\le i\le k$ and $\pi^+(i)=i$ for $k+1\le i\le n$.
In other words, the graph layout corresponding to $\pi^+$ starts with the vertices in $J$ arranged
according to $\pi$, followed by the vertices not in $J$ arranged in strictly increasing order.
Then the objective value of the FAS-QAP for $\pi^+$ consists of three parts:
the weight $W^{\pi}_1$ of the backward arcs going from $J$ into $J$, 
the weight $W_2$ of the backward arcs going from the complement of $J$ into the complement of $J$, and
the weight $W_3$ of the backward arcs going from the complement of $J$ into $J$.
We stress that the weights $W_2$ and $W_3$ only depend on $J$ but do not depend on the choice of $\pi$.
Hence we get for every permutation $\pi\in S_k$ that
\begin{equation}
\label{eq:principal.1}
\qap(A,F_n,\pi^+) ~=~ \qap(A[J],F_k,\pi)+W_2+W_3.
\end{equation}
Let $C$ be the $n\times n$ matrix in the linearization of the FAS-QAP for $A$.
Then
\begin{equation}
\label{eq:principal.2}
\lap(C,\pi^+) ~=~ \lap(C[J],\pi)+\sum_{i=k+1}^nc_{ii}.
\end{equation}
Equations (\ref{eq:principal.1}) and (\ref{eq:principal.2}) show that the FAS-QAP for $A[J]$ 
is linearizable.
The corresponding linearization matrix is $C[J]$ plus another linearization matrix that yields
a constant LAP objective value of $\sum_{i=k+1}^nc_{ii}-(W_2+W_3)$.
\qed

\begin{table}[tb]
\begin{center}
\begin{tabular}{|l||c|@{\quad}c@{\quad}|}
\hline && \\[-1.9ex]
 &$\qap(A[J],F_3,\pi)$
 &$\lap(C[J],\pi)$
\\[0.3ex]
\hline
&& \\[-1.9ex]
$\pi_1=(i,j,k)$ &$a_{ji}+a_{kj}+a_{ki}$ &$c_{ii}+c_{jj}+c_{kk}$\\[0.5ex]
$\pi_2=(i,k,j)$ &$a_{ji}+a_{jk}+a_{ki}$ &$c_{ii}+c_{jk}+c_{kj}$\\[0.5ex]
$\pi_3=(j,k,i)$ &$a_{ij}+a_{kj}+a_{ik}$ &$c_{ik}+c_{ji}+c_{kj}$\\[0.5ex]
$\pi_4=(j,i,k)$ &$a_{ij}+a_{kj}+a_{ki}$ &$c_{ij}+c_{ji}+c_{kk}$\\[0.5ex]
$\pi_5=(k,i,j)$ &$a_{ji}+a_{jk}+a_{ik}$ &$c_{ij}+c_{jk}+c_{ki}$\\[0.5ex]
$\pi_6=(k,j,i)$ &$a_{ij}+a_{jk}+a_{ik}$ &$c_{ik}+c_{jj}+c_{ki}$\\[1.5ex]
\hline
\end{tabular}
\end{center}
\caption{The objective values of the six permutations in the proof of Theorem~\protect{\ref{th:main.fas}}.}
\label{tab:six-permutations}
\end{table}

\begin{theorem}
\label{th:main.fas}
The FAS-QAP for matrix $A$ is linearizable, if and only if $A$ is a balanced 3-cycle matrix.
\end{theorem}
\proof
For the if part, we first use Theorem~\ref{th:subspace} to decompose $A$ into the sum 
of a symmetric matrix and a linear combination of directed cut matrices.
Lemmas~\ref{le:fas.1} and \ref{le:fas.2} imply that the FAS-QAP is linearizable for each of 
the summands, and then Lemma~\ref{le:linear} shows that the FAS-QAP is linearizable for 
matrix $A$ itself.  

For the only if part, consider a matrix $A$ for which the FAS-QAP is linearizable.
Lemma~\ref{le:fas.3} yields that the FAS-QAP for every principal $3\times3$ submatrix $A[J]$
defined by some $J=\{i,j,k\}$ with $i<j<k$ is linearizable.
We denote the corresponding linearization by $C[J]$, and for convenience we index the rows 
and columns of $C[J]$ also by $i<j<k$.
Table~\ref{tab:six-permutations} lists the objective values of the QAP and the LAP for 
the six permutations 
$\pi_1=(i,j,k)$, 
$\pi_2=(i,k,j)$,
$\pi_3=(j,k,i)$,
$\pi_4=(j,i,k)$,
$\pi_5=(k,i,j)$, and
$\pi_6=(k,j,i)$. 
Note that the sum of the LAP objective values for the three permutations $\pi_1$, $\pi_3$, $\pi_5$
equals the sum of LAP objective values for the three permutations $\pi_2$, $\pi_4$, $\pi_6$ 
(as both sums coincide with the sum of all the entries in matrix $C[J]$).
Consequently the two corresponding sums of QAP objective values are equal to each other as well,
which yields
\begin{eqnarray}
\lefteqn{(a_{ji}+a_{kj}+a_{ki}) \,+\, (a_{ij}+a_{kj}+a_{ik}) \,+\, (a_{ji}+a_{jk}+a_{ik}) ~=~}
\nonumber\\[0.5ex]
&=& (a_{ji}+a_{jk}+a_{ki}) \,+\, (a_{ij}+a_{kj}+a_{ki}) \,+\, (a_{ij}+a_{jk}+a_{ik})
\label{eq:fas-final}
\end{eqnarray}
Some algebraic simplifications turn (\ref{eq:fas-final}) into (\ref{eq:b3c}). 
As the choice of $i,j,k$ was arbitrary, matrix $A$ indeed is a balanced 3-cycle matrix.
\qed

\medskip
\section{The traveling salesman QAP}
\label{sec:tsp-qap}
An instance of the traveling salesman problem (TSP) consists of $n$ cities together with
an $n\times n$ distance matrix $A=(a_{ij})$.
The goal is to find a \emph{cyclic} permutation $\pi\in S_n$ that minimizes the linear
assignment function $\lap(A,\pi)$ in (\ref{eq:lap}).
We refer the reader to the book \cite{TSP} for a wealth of information on the TSP,
and to Burkard \& al \cite{BDDVW1998} for a survey on its well-solvable special cases.
The TSP can easily be formulated as a special case of the QAP, by choosing the first matrix $A$
as the underlying distance matrix and by choosing the second matrix as the $n\times n$ adjacency 
matrix $H_n=(h_{ij})$ of a directed Hamiltonian cycle whose entries are defined as follows: 
\begin{equation}
\label{eq:tsp}
h_{ij} ~=~ \left\{
\begin{array}{cl}
1 &\mbox{\quad if~ $j=i+1$, or if $i=n$ and $j=1$} \\[0.5ex]
0 &\mbox{\quad otherwise}
\end{array} \right.
\end{equation}
We will call this problem the \emph{TSP-QAP for matrix~$A$}, or just \emph{TSP-QAP} for short.
We stress that in the QAP formulation, all permutations $\pi\in S_n$ (and not just 
the cyclic ones) constitute feasible solutions.

In Theorem~\ref{th:main.tsp}, we will concisely characterize all linearizable instances of the TSP-QAP.
The proof of this theorem is based on the following result. 
\begin{proposition}
\label{pr:const-tsp}
(Gabovich \cite{Gabovich1976}, and independently Berenguer \cite{Be1979})~
The following two statements are equivalent:
\begin{itemize}
\itemsep=0.0ex
\item[(i)]  For the distance matrix $A$, all permutations $\pi$ yield the same TSP objective value.
\item[(ii)] Matrix $A$ is a weak sum matrix.  \qed
\end{itemize}
\end{proposition}
Gilmore, Lawler \& Shmoys \cite{GLS1985} present a very simple and concise proof of 
Proposition~\ref{pr:const-tsp} by means of linear algebra.

\begin{theorem}
\label{th:main.tsp}
The TSP-QAP for matrix $A$ is linearizable, if and only if $A$ is a weak sum matrix.
\end{theorem}
\proof
For the if part, we assume without loss of generality that $A$ is a sum matrix.
Then by Proposition~\ref{pr:const-tsp} all permutations yield the same QAP objective value,
and it can be linearized trivially by a matrix $C$ that yields the same constant LAP 
objective value for all permutations.  

For the only if part, consider an $n\times n$ matrix $A$ for which the TSP-QAP is 
linearizable and let $C$ be the corresponding linearization.
For a permutation $\pi\in S_n$, its cyclic shift is the permutation $\pi\xxx{1}$ defined by 
$\pi\xxx{1}(i)=\pi(i+1)$ for $1\le i\le n-1$ and $\pi\xxx{1}(n)=\pi(1)$.
For $0\le k\le n-1$, the $k$th cyclic shift of $\pi$ results by cyclically shifting it $k$ times;
note that $\pi\xxx{0}=\pi$.
Now let us consider the total objective value of all $n$ cyclic shifts $\pi\xxx{0},\ldots,\pi\xxx{n-1}$
of permutation $\pi$ for QAP and LAP.
In the QAP, every cyclic shift $\pi\xxx{k}$ has the same objective value.
All cyclic shifts correspond to the same tour through the cities, and they only differ in 
the choice of their starting point.
This yields
\begin{equation}
\label{eq:tsp.1}
\sum_{k=0}^{n-1} \qap(A,H_n,\pi\xxx{k}) ~=~ n\cdot \qap(A,H_n,\pi).
\end{equation}
In the LAP, the $n$ shifts cover every element of matrix $C$ exactly once.
This yields
\begin{equation}
\label{eq:tsp.2}
\sum_{k=0}^{n-1} \lap(C,\pi\xxx{k}) ~=~ \sum_{i=1}^n \sum_{j=1}^n c_{ij}.
\end{equation}
Since the values in (\ref{eq:tsp.1}) and (\ref{eq:tsp.2}) coincide, this implies that all
tours in the traveling salesman have the same length $(\sum_{i=1}^n \sum_{j=1}^n c_{ij})/n$.
Then Proposition~\ref{pr:const-tsp} yields that $A$ indeed is a weak sum matrix.
\qed

\medskip
\section{Conclusion}
\label{sec:conclusion}
We have given combinatorial characterizations of the linearizable instances for two classes 
of asymmetric QAPs: the weighted feedback arc set QAP, and the traveling salesman QAP.
Similarly as in the symmetric case, all these linearizable asymmetric instances carry a 
very strong and very restrictive combinatorial structure.

Our results (together with the known results on the symmetric case) might indicate that
linearizable instances of the QAP are rare events and will essentially never show up in 
real world situations.
It would be interesting to support these intuitions by means of a probabilistic analysis
in some reasonable stochastic model.

Another line for future research is to identify further linearizable families for the
asymmetric case.
A more ambitious goal would be to get a complete combinatorial characterization of
all linearizable asymmetric QAP instances.

\bigskip
\small
\paragraph{Acknowledgements.}
Part of this research was conducted while Vladimir Deineko and Gerhard Woeginger were visiting TU Graz,
and they both thank the Austrian Science Fund (FWF): W1230, Doctoral Program in ``Discrete Mathematics'' 
for the financial support.
Vladimir Deineko acknowledges support
by Warwick University's Centre for Discrete Mathematics and Its Applications (DIMAP).
Gerhard Woeginger acknowledges support
by DIAMANT (a mathematics cluster of the Netherlands Organization for Scientific Research NWO).
\normalsize


\end{document}